\documentclass[12pt]{amsart}
\usepackage{amssymb, amstext, amscd, amsmath}
\makeatletter
\def\@cite#1#2{{\m@th\upshape\bfseries%
[{#1\if@tempswa{\m@th\upshape\mdseries, #2}\fi}]}}
\makeatother
\theoremstyle{plain}
\newtheorem{thm}{Theorem}[section]

\newtheorem{prop}[thm]{Proposition}

\theoremstyle{definition}
\newtheorem{rem}[thm]{Remark}
\newtheorem{defn}[thm]{Definition}

\newtheorem{eg}[thm]{Example}
\newcommand{\Prf}{\noindent\textbf{Proof.\ }}
\newcommand{\bx}{\hfill$\blacksquare$\medbreak}

\newcommand{\ca}{\mathrm{C}^*}

\newcommand{\bbC}{{\mathbb{C}}}

  \newcommand{\A}{{\mathcal{A}}}
  \newcommand{\B}{{\mathcal{B}}}

  \newcommand{\F}{{\mathcal{F}}}
  
\renewcommand{\H}{{\mathcal{H}}}
  \newcommand{\I}{{\mathcal{I}}}

  \newcommand{\M}{{\mathcal{M}}}

\newcommand{\fA}{{\mathfrak{A}}}

\newcommand{\qwhere}{\quad\text{where}\quad}

\newcommand{\dom}{\operatorname{Dom}}

\newcommand{\ran}{\operatorname{Ran}}
\renewcommand{\dom}{\operatorname{Dom}}

\newcommand{\trace}{\operatorname{trace}}

\def\one{{\mathchoice{\rm 1\mskip-4mu l}{\rm 1\mskip-4mu l}{\rm 1\mskip-4.5mu l}{\rm
1\mskip-5mu l}}}

\begin{document}

\title[]{A unified framework for graph algebras and quantum causal histories}
%
%
\author[D.W. Kribs]{David~W.~Kribs}
\address{Department of Mathematics and Statistics, University of
Guelph, Guelph, Ontario, Canada  N1G 2W1} \address{Institute for
Quantum Computing, University of Waterloo, Waterloo, ON, CANADA
N2L 3G1} \address{Perimeter Institute for Theoretical Physics, 31
Caroline St. N., Waterloo, ON, CANADA N2L 2Y5}
\begin{abstract}
We present a mathematical framework that unifies the quantum
causal history formalism from theoretical high energy physics and
the directed graph operator framework from the theory of operator
algebras. The approach involves completely positive maps and
directed graphs and leads naturally to a new class of operator
algebras.
\end{abstract}

\thanks{2000 {\it  Mathematics Subject Classification.} 46L05, 47L75, 81P68, 83C45.}
\thanks{{\it Key words and phrases.} directed graph, completely positive
map, partial isometry, quantum operation,  quantum causal
history.}
\maketitle

\section{Introduction}   \label{S:intro}

In this paper we present a new operator theoretic framework that
provides a unified approach for recent studies in theoretical high
energy physics and contemporary operator algebra theory. More
specifically, this approach involves completely positive maps and
directed graphs, and includes the quantum causal history formalism
from recent work towards a quantum theory of gravity on the one
hand and the graph-operator framework from the theory of directed
graph operator algebras on the other. We also define a new class
of operator algebras that is naturally defined through this
approach.

The basic physical properties that a quantum theory of gravity
must satisfy motivated F. Markopoulou to invent a formalism called
``Quantum Causal Histories'' \cite{HMS,M02,M00C,M00B,M00A}. A
secondary goal of this paper is to introduce this formalism to the
broader mathematical community. Fundamental examples include
causal evolutions of spin networks \cite{Gupta,MS97} and quantum
computers \cite{M00C}. The basic definitions have been refined
through a series of papers and now a clean mathematical definition
is emerging \cite{HMS}. Mathematically, and somewhat roughly
speaking at this point, a {\it quantum causal history} (QCH) is
given by a directed graph with a finite-dimensional Hilbert space
for each vertex and a quantum operation associated with each
directed edge. The vertex spaces correspond to events, or
observables, within a local history and the quantum operations
indicate causal links between pairs of related events. As
described below, the QCH framework incorporates further structure
motivated by the characterization of evolution within quantum
systems.

On the other hand, Cuntz and Krieger \cite{CK} were motivated by a
problem in dynamical systems to initiate the study of operator
algebras associated with directed graphs. Over the past two
decades the study of directed graph operator algebras and related
topics has taken on a life of its own and now, it is fair to say,
plays a central role in the modern theory of operator algebras. We
mention the articles \cite{BP} -- \cite{Spi} as entrance points
into the extensive literature on the subject.

In $\S~2$ we recall some basic properties of completely positive
maps and quantum operations. We define the ``CP -- directed graph
framework'' and associated $\ca$-algebras in $\S~3$, and show how
graph algebras fit into the framework. In $\S~4$ we discuss in
some detail the QCH formalism, draw a connection with quantum
computing, and prove a theorem for the QCH $\ca$-algebras.

\section{Completely Positive Maps}   \label{S:cpmaps}

Given a Hilbert space $\H$ we let $\B(\H)$ be the set of bounded
operators that act on $\H$. A {\it completely positive} (CP) map
is a linear map $\Phi : \B(\H_1)\rightarrow\B(\H_2)$ such that the
``ampliated'' maps
\[
 \one_k \otimes \Phi : \M_k \otimes \B(\H_1) \longrightarrow
\M_k \otimes \B(\H_2)
\]
are positive for $k\geq 1$. (Here $\M_k$ denotes the set of
$k\times k$ complex matrices and $\one_k$ denotes the identity
operator, the ``maximally mixed state'', inside $\M_k$.) For basic
mathematical properties of CP maps see \cite{Paulsentext2} and
physical motivations see \cite{Kraustext}.

A fundamental technical device in the study of CP maps is the {\it
operator-sum representation} given by the theorem of Choi
\cite{Choi} and Kraus \cite{Kraus}. For every CP map $\Phi$ on
finite-dimensional space, there is a set of {\it noise operators}
(or {\it errors}) $\{A_i\}\subseteq \B(\H_1,\H_2)$ such that
\begin{eqnarray}\label{opsum}
\Phi(\rho) &=& \sum_i A_i\, \rho\, A_i^* \quad\forall
\rho\in\B(\H_1).
\end{eqnarray}
The map is {\it unital} if also $\sum_i A_i A_i^* = \one_{\H_2}$.
It is always possible to choose a family of noise operators with
cardinality at most $\dim(\H_1)\dim(\H_2)$. On infinite
dimensional space not all CP maps have such a form, and when they
do the sum in (\ref{opsum}) converges in the strong operator
topology. For brevity, we assume all the CP maps considered here
have a representation as in (\ref{opsum}).

A {\it quantum operation} (or {\it quantum evolution}, or {\it
quantum channel}) is a CP map $\Phi : \B(\H_1)\rightarrow\B(\H_2)$
that is also trace preserving. When $\Phi$ is represented as in
(\ref{opsum}), trace preservation is equivalent to the identity
\begin{eqnarray}\label{tracepreserve}
 \sum_i A_i^*  A_i = \one_{\H_1}.
\end{eqnarray}
Thus, a quantum operation $\Phi$ is a map that satisfies
(\ref{opsum}) and (\ref{tracepreserve}) for some set of operators
$\{A_i\}$. Equivalently, a quantum operation is a CP map such that
its associated {\it dual map}, denoted by $\Phi^\dagger:
\B(\H_2)\rightarrow\B(\H_1)$, is unital. (Recall that the dual map
for a map $\Phi$ is defined via the equation
$\trace(\rho\,\Phi^\dagger(\sigma)) =
\trace(\Phi(\rho)\,\sigma)$.)

The ideal physical examples of quantum operations are {\it unitary
maps}  as they characterize evolution of states within a closed
quantum system. Such a map is of the form $\Phi(\rho) = U \rho\,
U^*$ for some unitary operator $U$. When evolution occurs in an
open system (i.e., when the system of interest is exposed to an
external environment) quantum operations have the more general
form given by (\ref{opsum}) and (\ref{tracepreserve}). See
\cite{Ksurvey} for further discussions and references.

\section{The CP -- Directed Graph Framework}   \label{S:cpgraph}

Let $E=(E^0,E^1,r,s)$ be a (countable) directed graph with
vertices $x\in E^0$, directed edges $e\in E^1$ and range and
source maps $r,s: E^1\rightarrow E^0$ giving the initial ($s(e)$)
and final ($r(e)$) vertices of an edge $e$.  When $e\in E^1$
satisfies $s(e)=x$ and $r(e)=y$, we shall write $e=(x,y)$.

Suppose we have a Hilbert space $\{\H(x):x\in E^0\}$ for each
vertex and a family of CP maps $\Psi=\{\Phi_e:e\in E^1\}$ with
domains and ranges that satisfy
\[
(\dagger)\left\{
\begin{array}{lll}
(i) & \dom(\Phi_e) = \B(\H(s(e))) & \mbox{$\forall\, e \in E^1$} \\
(ii) & \ran(\Phi_e) \subseteq \B(\H(r(e))) & \mbox{$\forall\, e\in
E^1$}
\end{array}\right.
\]

Given such a family of spaces and maps, define the Hilbert space
$\H = \oplus_{x\in E^0} \H(x)$ and let $P_x$ be the projection of
$\H$ onto $\H(x)$.

\begin{defn}
Given a directed graph $E$, let $\{\H(x):x\in E^0\}$ be Hilbert
spaces and let $\Psi=\{\Phi_e: e\in E^1\}$ be a family of CP maps
that satisfy $(\dagger)$. Suppose $\Phi_e = \{A_{e,i}:i\in\I_e\}$
is an operator-sum representation of $\Phi_e$ for each $e\in E^1$.
We can naturally regard the $A_{e,i}$ as operators acting on $\H$.
Define $\fA_\Psi$ to be the $\ca$-algebra generated by all
operators $A_{e,i}$ and vertex projections $P_x$; so that,
\[
\fA_\Psi = \ca \big( \{ P_x, A_{e,i} : x\in E^0,\, e\in E^1, \,
i\in\I_e\}\big).
\]
\end{defn}

The choice of noise operators that represent a given CP map in
(\ref{opsum}) is of course not unique. However, as our notation
suggests the algebras $\fA_\Psi$ are independent of these choices.

\begin{prop}
Let $E$ be a directed graph and let $\Psi= \{\Phi_e:e\in E^1\}$ be
a family of CP maps that satisfy $(\dagger)$. Then the algebra
$\fA_\Psi$ is independent of the choice of operators $\{A_{e,i}\}$
that represent the maps $\Phi_e$ as in (\ref{opsum}).
\end{prop}

\Prf Suppose that $\{A_{e,i}\}_i$ and $\{A_{e,j}^\prime\}_j$
represent $\Phi_e$ via  equation (\ref{opsum}). By possibly
including zero operators we may assume the cardinality of these
two sets is the same. Then from the structure theory for CP maps,
there is a scalar unitary matrix $U=(u_{ij})$ such that
\[
A_{e,i} = \sum_j u_{ij} A_{e,j}^\prime \quad \forall i.
\]
It follows that the algebras generated by
$\{P_x,A_{e,i}\}_{x,e,i}$ and $\{P_x,A_{e,j}^\prime\}_{x,e,j}$
coincide, and the result follows.
 \bx

\subsection{Graph Algebras}\label{sS:graphalg}

We now discuss one of the motivating special cases for this
framework. Let $E=(E^0,E^1,r,s)$ be a directed graph. Consider
families of operators $\{P_x,S_e: x\in E^0, e\in E^1\}$, where the
$P_x$ are projections and the $S_e$ are partial isometries (or
equivalently, unitary operators restricted to a subspace), that
act on the same Hilbert space and satisfy:
\[
(\ddagger)\left\{
\begin{array}{lll}
(i) & S_{e}^{*}S_e = P_{s(e)} & \mbox{$\forall\, e \in E^1$} \\
(ii) & S_{e}S_e^* \leq P_{r(e)} & \mbox{$\forall\, e\in E^1$}
\end{array}\right.
\]
Then the structure of $E$ determines the relations satisfied by
$\{P_x,S_e\}$ in the sense that the initial projection for each
$S_e$ is equal to the projection for the source vertex of $e$ and
the range projection for each $S_e$ is supported on the projection
for the range vertex of $e$.

The relations $(\ddagger)$ provide the fundamental base case for
investigations into operator algebras associated with directed
graphs. In the most general context, a {\it graph algebra} is an
operator algebra generated by a family $\{S_e,P_x\}$. There are a
number of refinements and generalizations of the formulation
$(\ddagger)$. In many instances the $S_e$ are assumed to have
mutually orthogonal ranges. The projections $P_x$ are typically
assumed to have mutually orthogonal ranges as well, or sometimes
just mutually commuting ranges. There are also topological graph
generalizations wherein the vertices and edges are locally compact
spaces and the range and source maps are continuous maps. However,
in every setting the motivating case is the same: A Hilbert space
$\H(x)$ associated with every vertex $x$ in $E$ and for every
directed edge $e=(x,y)$ a partial isometry $S_e$ that maps from
$\H(x)$ to $\H(y)$.

If $\{S_e,P_x\}$ satisfy $(\ddagger)$, observe that for each
$e=(x,y)$ the operator $S_e$ defines a unitary from
$\H(x)=P_{x}\H$ into $\H(y)=P_{y}\H$ and a unitary CP map
$\Phi_e:\B(\H(x))\rightarrow\B(\H(y))$ via
\[
\Phi_e(\rho) = S_e^\prime\, \rho\, (S_e^\prime)^* \qwhere
S_e^\prime = S_e|_{\H(x)}.
\]
Thus the corresponding algebra $\fA_\Psi$ defines a graph algebra,
and so graph algebras form a subclass of the algebras $\fA_\Psi$.

\section{The Quantum Causal History Formalism}   \label{S:qch}

The mathematical formalism for QCH's has undergone a series of
refinements since being introduced in \cite{M00B}. The
presentation below is most closely related to the recent
formulation of Hawkins, Markopoulou and Sahlmann \cite{HMS}. The
nomenclature we use is slightly different than \cite{HMS}, we do
this to mesh with the graph algebra terminology. We shall focus on
the mathematical aspects and touch on the physical motivations for
various constraints.

To define a QCH then, we begin with a graph $E=(E^0,E^1,r,s)$,
which may also be interpreted as a partial order when there are no
loops. This graph represents a causal set wherein the vertices
correspond to a set of local events in the universe and vertices
linked by directed edges indicate causal relations between events.
From the postulates of quantum mechanics, events are represented
by density operators on Hilbert space. Recent work in string
theory and loop quantum gravity (see \cite{M02}) suggests that any
finite region of space should contain a finite amount of
information. Thus, each of the event spaces is assumed to be
finite-dimensional. Since causality can be interpreted as
transferring information from one event to another and because, by
definition, a QCH describes local causality at the quantum level,
a causal relation given by a directed edge $e=(x,y)\in E^1$
corresponds to a quantum operation $\Phi(x,y) : \B(\H(x))
\rightarrow \B(\H(y))$ between event spaces.

Thus, at its mathematical core, a QCH consists of a directed
graph, with a finite-dimensional Hilbert space for each vertex,
and a quantum operation for each directed edge. There are further
constraints within a QCH and we discuss them now briefly.

First some terminology. Given $x,y\in E^0$ write $x\leq y$ when
$x$ precedes $y$ as an event. In this case there is a
future-directed curve from $x$ to $y$ and this is represented by a
directed edge $e=(x,y)\in E^1$. If $x\leq y$ or $y\leq x$ then $x$
and $y$ are {\it related} and otherwise they are {\it spacelike
separated} and we use $x\sim y$ to denote this.  A path in $E$
corresponds to a {\it future-directed path} through the events in
the history. Such a path is {\it future (past) inextendible} if
there is no event in $E$ which is in the future (past) of the
entire path. Loops in $E$ correspond to {\it closed timelike
curves}. From the finiteness assumption discussed above, $E$ is
{\it locally finite} in the sense that for any $x,y\in E^0$ there
are at most finitely many $z\in E^0$ such that $x\leq z \leq y$.
Given $x,y\in E^0$, there is also  no generality lost in assuming
there is at most one edge $e=(x,y)$ in $E$ from $x$ to $y$. (If
$s(e)=x=s(f)$ and $r(e)=y=r(f)$ then the operations associated
with these edges could be combined to form a single operation that
encodes the relevant causal structure from event $x$ to event
$y$.)

An {\it acausal set} $\xi\subseteq E^0$ is defined by the property
that $x\sim y$ whenever $x,y\in\xi$. Such a set is a {\it complete
future} for an event $x$ if $\xi$ intersects any future
inextendible future-directed path that starts at $x$. A {\it
complete past} is defined analogously. The composite state space
for $x\sim y$ (the physical existence of which is guaranteed by
quantum mechanics) is $\H(\{x,y\}) = \H(x) \otimes \H(y)$ and more
generally $\H(\xi) = \otimes_{x\in\xi} \H(x)$. For $x\in E^0$
write $\A(x)$ for the matrix algebra $\B(\H(x))$ and similarly
define $\A(\xi)=\otimes_{x\in\xi}\A(x)$ for a set $\xi\subseteq
E^0$. Given an acausal set $\xi$ and an event $x\in \xi$, there is
a natural unital embedding $\iota_x: \A(x)\hookrightarrow\A(\xi)$
and we shall write $\A(x)\subseteq\A(\xi)$.

If $\xi$ and $\zeta$ are acausal sets such that $\xi$ is a
complete past for $\zeta$ and $\zeta$ is a complete future for
$\xi$, then we write $\xi \preceq \zeta$ and say that
$(\xi,\zeta)$ form a {\it complete pair}. Such a pair represents
an evolution in a closed quantum system, hence woven into the
fabric of the QCH there should be a unitary operator
$U(\xi,\zeta): \H(\xi) \rightarrow \H(\zeta)$. Such an operator
determines a unitary map (an isomorphism) $\Phi(\xi,\zeta):
\A(\xi)\rightarrow\A(\zeta)$ via
\[
\Phi(\xi,\zeta)(\rho) = U(\xi,\zeta)\, \rho\, U(\xi,\zeta)^*
\quad\forall \,\rho\in\A(\xi).
\]
Note the restriction of $\Phi(\xi,\zeta)$ (respectively
$\Phi(\xi,\zeta)^\dagger$) to $\A(x)\subseteq\A(\xi)$ for
$x\in\xi$ (respectively $\A(z)\subseteq\A(\zeta)$ for $z\in\zeta$)
is a $\ast$-homomorphism. This gives the structure of a QCH at the
global level, but does not indicate how the isomorphisms
$\Phi(\xi,\zeta)$ should depend on the individual causal relations
between events in $\xi$ and $\zeta$. This is the role played by
the operations $\Phi(x,y)$ on individual edges.

We now give a precise mathematical definition of a QCH. We note
that the maps in the definition below have directions reversed
from the presentation in \cite{HMS}. Here we take the dual
approach so the ``directions'' of the maps are in line with the
graph structure. Recall that if $\iota_A: \A_1\rightarrow\A_2$ and
$\iota_B: \B_1\rightarrow\B_2$ are $\ast$-monomorphisms and $\Psi:
\A_2\rightarrow\B_2$ is a map, then the {\it reduction} of $\Psi$
to $\A_1\mapsto\B_1$ is the map $\Phi = \iota_B^\dagger \circ \Psi
\circ i_A$.

\begin{defn}
A {\it quantum causal history} consists of a directed graph $E$
with a Hilbert space $\{\H(x): x\in E^0\}$ for each event and a
quantum operation $\Phi (x,y) : \A(x)\rightarrow\A(y)$ for each
pair of related events $x\leq y$ such that the following axioms
are satisfied:

$(i)$ (Extension) For all $y\in E^0$ and $\zeta\subseteq E^0$ a
complete future of $y$, there is a homomorphism $\Phi_F(y,\zeta) :
\A(y)\rightarrow\A(\zeta)$ such that $\Phi_F(y,\zeta)^\dagger$ is
a quantum operation and for all $z\in\zeta$, the reduction of
$\Phi_F(y,\zeta)^\dagger$ to $\A(z)$ is $\Phi(y,z)^\dagger$.
Likewise, for all $y\in E^0$ and $\xi\subseteq E^0$ a complete
past of $y$, there is a quantum operation
$\Phi_P(\xi,y):\A(\xi)\rightarrow\A(y)$ such that
$\Phi_P(\xi,y)^\dagger$ is a homomorphism and for all $x\in\xi$
the reduction of $\Phi_P(\xi,y)^\dagger$ to $\A(x)$ is
$\Phi(x,y)^\dagger$.

$(ii)$ (Spacelike Commutativity)  If $x\sim y\in E^0$ and
$\zeta\subseteq E^0$ is a complete future of $x$ and $y$, then the
images of $\Phi_F(x,\zeta)$ and $\Phi_F(y,\zeta)$ commute inside
$\A(\zeta)$.  Likewise, if $y\sim z \in E^0$ and $\xi\subseteq
E^0$ is a complete past of $y$ and $z$, then the images of
$\Phi_P(\xi,y)^\dagger$ and $\Phi_P(\xi,z)^\dagger$ commute inside
$\A(\xi)$.

$(iii)$ (Composition) If $\zeta\subseteq E^0$ is a complete future
of $x$ and a complete past of $y$, then $\Phi(x,y) =
\Phi_P(\zeta,y)\circ\Phi_F(x,\zeta)$.
\end{defn}


If $\xi\preceq\zeta$ form a complete pair within a given QCH, it
is proved in \cite{HMS} that there is a unique unitary map
$\Phi(\xi,\zeta):\A(\xi)\rightarrow\A(\zeta)$ such that the
reduction of $\Phi(\xi,\zeta)$ to $\A(x)\mapsto\A(y)$  is
$\Phi(x,y)$ for all $x\in\xi$ and $y\in\zeta$. Thus the
isomorphisms $\Phi(\xi,\zeta)$ discussed above may be built up
from the individual edge maps $\Phi(x,y)$ and hence the edge maps
are the fundamental building blocks for a QCH.

\begin{eg}\label{S:qcomp}

({\bf Quantum Computers}) As discussed in \cite{M00C}, the basic
model for a quantum computer fits into the QCH formalism.
Specifically, each quantum algorithm may be interpreted as a QCH
via its ``circuit-gate'' presentation. (See \cite{Ksurvey} for a
brief mathematical introduction to quantum algorithms and
references.) The QCH for a given algorithm has vertex spaces all
equal to $\bbC^2$. The directed edges correspond to the choice of
unitary gates within the algorithm and the vertex spaces encode
the intermediate states of the quantum bits of information (the
`qubits'). The structure of the associated directed graph is the
same as the circuit-gate diagram, with the circuits labelled as
vertices and the gates labelled as directed edges.
\end{eg}

We finish by proving that the algebras $\fA_\Psi$ associated with
QCH's are familiar objects from operator theory. For basic
properties of AF-algebras we point the reader to the text
\cite{byeg}.

\begin{thm}\label{qchopalg}
Let $\fA_\Psi$ be the $\ca$-algebra associated with a given
quantum causal history. Then $\fA_\Psi$ is an AF-algebra.
\end{thm}

\Prf It is enough to prove that every finite set of elements of
$\fA_\Psi$ can be approximated by elements lying in a
finite-dimensional subalgebra. But elements of the form $A =
A_1\cdots A_n$, where each $A_k=A_{e,i}$ or $A_k = A_{e,i}^*$ for
some $e$ and $i$, span a dense subspace of $\fA_\Psi$. (Note that
the vertex projections are obtained via equation
(\ref{tracepreserve}).) Hence it is enough to show that each
finite set of such elements lies in a finite-dimensional
subalgebra. Suppose $\F$ is such a set. Note that each $A_{e,i}$
belongs to $\B(\H(x),\H(y))$ for some $x,y\in E^0$, and so the
same is true for every $A_{e,i}^*$ and all elements $A\in\F$.
Thus, as $\dim\H(x)<\infty$ for all $x\in E^0$, it follows that
the algebra generated by $\F$ is a finite-dimensional algebra
which is also a subalgebra of $\fA_\Psi$, and the result follows.
\bx

\begin{rem}
While graph algebras and quantum causal histories provided the
initial impetus for the CP -- directed graph framework presented
here, it is evident that this structure admits other
possibilities. Indeed, there are many CP maps that are neither
unitary maps nor quantum operations, and presumably the
$\ca$-algebras $\fA_\Psi$ would go beyond the graph algebra and
QCH subclasses in general. We have also not considered here the
various possibilities for non-selfadjoint algebras defined by this
framework. We plan to undertake investigations of these algebras
elsewhere and hope this paper motivates others to do the same.
\end{rem}



\vspace{0.05in}

{\noindent}{\it Acknowledgements.}
I would like to thank Fotini Markopoulou
for enlightening conversations and Eli Hawkins for helpful
comments on an early draft. I am also grateful to colleagues at
the Institute for Quantum Computing and the Perimeter Institute
for interesting discussions. This work was partially supported by
an NSERC grant.

\end{document}